\documentclass[11pt]{amsart}
\usepackage{amsfonts,amssymb,amsmath}

\newcommand{\Rm}{\mathbb{R}}

\newcommand{\mL}{\mathcal{L}}

\newcommand{\Vm}{\ensuremath{\mathbb{V}}}
\newcommand{\mM}{\ensuremath{\mathcal{M}}}

\newtheorem{lem}{Lemma}
\newtheorem{thm}[lem]{Theorem}

\newtheorem{prop}[lem]{Proposition}

\def\qed {\mbox{}\hfill {\small \fbox{}} \\}
\def\lto{\longrightarrow}
\def\lmto{\longmapsto}

\def\leq{\leqslant}

\title[Flows]
{Uniqueness of signed measures solving \\ the continuity equation
for \\ Osgood vector fields}

\author{Luigi Ambrosio, Patrick Bernard}
\address{
{\rm Luigi Ambrosio}\\
Scuola Normale Superiore\\
Piazza dei Cavalieri 7\\
56123 Pisa, Italy}
\email{l.ambrosio@sns.it}
\address{
{\rm Patrick Bernard}\\
Universit\'e Paris-Dauphine\\
et CEREMADE, UMR CNRS 7534\\
Pl. du Mar\'echal de Lattre de Tassigny\\
75775 Paris Cedex 16\\
France (Membre de l'IUF)} 
\email{Patrick.Bernard@ceremade.dauphine.fr}
\date{}
\begin{document}
\begin{abstract}
Nonnegative measure-valued solutions of the continuity equation are
uniquely determined by their initial condition, if the
characteristic ODE associated to the velocity field has a unique
solution. In this paper we give a partial extension of this result
to signed measure-valued solutions, under a quantitative two-sided Osgood
condition on the velocity field. Our results extend those obtained
for log-Lipschitz vector fields in \cite{BC}.
\end{abstract}

\maketitle

\section{Introduction}

Let $T>0$ and let
$$
V(t,x):(0,T)\times \Rm^d\lto \Rm^d
$$
be a Borel vectorfield.
We associate to $V$ the equations
\begin{equation}\tag{ODE}
\dot \gamma(t)=V(t,\gamma(t))
\end{equation}
and (with the notation $V_t(x)=V(t,x)$)
\begin{equation}\tag{PDE}
\partial_t \mu_t + \text div(V_t \mu_t)=0.
\end{equation}
A solution of (ODE) is an absolutely continuous curve $\gamma(t)$
such that $\dot \gamma(t)=V(t,\gamma(t))$ almost everywhere on
$[0,T]$.
We shall also consider generalized solutions in the sense
 of Filippov, see more details below. The so-called \emph{continuity} or \emph{Liouville} equation (PDE)
is considered in the sense of distributions. We
shall work with solutions in the class of measures. We denote by
$\mM(\Rm^d)$ the set of signed Borel measures with finite total
variation on $\Rm^d$, by $\mM^+(\Rm^d)$ the subset of non-negative
finite measures, and by $|\mu|\in \mM^+(\Rm^d)$ the total variation
of a measure $\mu\in \mM(\Rm^d)$. We shall consider only solutions
$\mu_t$ satisfying $|\mu_t|(\Rm^d)\in L^\infty(0,T)$; this is not a
very restrictive assumption, because many approximation schemes do
provide solutions $\mu_t$ with this property. We shall also assume
that
\begin{equation}\tag{I}
 \int_0^T\int_{\Rm^d} \|V_t\| d|\mu_t|dt<\infty,
\end{equation}
a property surely satisfied if $\|V\|$ is uniformly bounded. Under
these assumptions the notion of distributional solution is well
defined, and it is equivalent to the requirement that, for all
$\phi\in C^1_c(\Rm^d)$, $t\mapsto\int\phi\,d\mu_t$ belongs to the
Sobolev space $W^{1,1}(0,T)$, with distributional derivative given
by
$$
\int_{\Rm^d}\langle V_t(x),\nabla\phi(x)\rangle\,d\mu_t(x).
$$
Using this fact, and the uniform continuity properties of Sobolev
functions on the real line, it is easy to check (see for instance
\cite[Proposition~8.1.7]{AGS}) that we can restrict ourselves
(possibly modifying $\mu_t$ in a negligible set of times) to weakly
continuous solutions $t\mapsto\mu_t$, in the duality with
$C_c(\Rm^d)$. Moreover, the initial condition $\mu_0$ for (PDE) is
defined in a weak sense:
$$
\lim_{t\downarrow
0}\int_{\Rm^d}\phi\,d\mu_t=\int_{\Rm^d}\phi\,d\mu_0\qquad \forall
\phi\in C_c(\Rm^d).
$$
So, from now on only weakly continuous solutions $\mu_t$ will be
considered. Reversing the time variable, also the final condition
$\mu_T$ is well defined, still in the weak sense.

Our goal in the present paper is to study the relations between
uniqueness for (ODE) and uniqueness for (PDE). It is known that
uniqueness for (ODE) implies, via the so-called superposition
principle, that nonnegative solutions of (PDE) are uniquely
determined by the initial condition $\mu_0$, see
\cite{AGS,Am:cetraro,AC,maniglia,B}. The question turns out to be
much more subtle if we work in the class of signed measures. Of
course, if $\mu_t$ is a solution, we can write it as the difference
of the two non-negative measures $\mu_t^+$ and $\mu_t^-$. However,
these measures need not solve the equation. This remark is
reminiscent of the notion of renormalized solutions, see
\cite{DPL,Am:04}: we may call renormalized a solution $\mu_t$ such
that $\mu_t^+$ and $\mu_t^-$ are both solutions (or equivalently
such that $|\mu_t|$ is a solution). It is clear that there is
uniqueness if all distributional solutions are renormalized, and if
there is uniqueness for (ODE), but
 renormalized solutions have been studied only under
weak differentiability assumptions on $V_t$, and only in the class
of absolutely continuous measures $\mu_t$ (see \cite{Am:cetraro} for
a survey on this topic).
In this paper, we leave aside 
the question of the general relations between (ODE) uniqueness
and (PDE) uniqueness, and we focus on a particular class
of vectorfields for which (ODE) uniqueness is well-known, and
derive some consequences at the (PDE) level (and, in particular, that
all solutions are renormalized).

We recall that a modulus of continuity is a continuous
non-decreasing function $\rho:[0,1)\lto [0,\infty)$, such that
$\rho(0)=0$. A modulus of continuity $\rho$ is said to be Osgood if
$$
\int_0^1 \frac{1}{\rho(s)} ds=+\infty.
$$
We will always extend the moduli of continuity to $[1,\infty)$ by
$\rho=\infty$. Typical examples of Osgood moduli of continuity are
$\rho(s)=s$ and $\rho(s)= s(1-\ln(s))$. Note that the moduli
$\rho(s)=s^\alpha$, $\alpha\in (0,1)$, are not Osgood.

It is known that uniqueness holds for (ODE) if there exist a Osgood
modulus of continuity $\rho$ and $C\in L^1(0,T)$ such that
\begin{equation}\tag{O}
|\langle V(t,x)-V(t,y),x-y\rangle|\leq
C(t)\|x-y\|\rho(\|x-y\|)
\end{equation}
for all $x,\,y\in\Rm^d$, and all $t\in (0,T)$. Condition (O) does not
seem to imply continuity of $V_t$: in the case
when the modulus $\rho$ is linear, (O) implies that the symmetric part
of the distributional derivative is bounded, hence Korn's inequality gives that
$V_t$ is equivalent, up to Lebesgue negligible sets, to a continuous function. Since 
we consider measures $\mu_t$ that are possibly singular, even in the case when $\rho$
is linear we can not apply this result to reduce ourselves to a continuous vector field;
therefore we will not investigate the continuity question here (also because
adding the continuity assumption would not lead to a great simplification of the
uniqueness proof).

In order to prevent blow-up of solutions, the following bound is useful:
\begin{equation}\tag{B}
| V(t,x)|\leq D(t)\quad\forall x\in\Rm^d,\,\,\forall t\in
(0,T),\,\,\text{for some $D\in L^1(0,T)$.}
\end{equation}

The equation (ODE) is well understood under (O) and (B):
There exists a unique  flow map 
$$X:[0,T]\times [0,T] \times \Rm^d \lmto \Rm^d$$
which is such that 
$X(s,t,\cdot)$ is a homeomorphism of $\Rm^d$ for each $s$ and $t$;
$X(t,t,\cdot)=Id$ for each $t$. In addition
$$t\mapsto X(s,t,x)$$
is  a generalized solution of (ODE) in the sense of Filippov
(the definition is recalled below)
 for each $s$ and $x$.
In the case when $V_t$ is continuous, then generalized 
solutions in the sense of Filippov are just ordinary solutions of (ODE).
Uniqueness
implies that $X$ satisfies the semigroup property
\begin{equation}\label{semigroup}
X(t_3,t_2,X(t_1,t_3,x))=X(t_1,t_2,x)\quad
\forall x\in\Rm^d,\,\,\forall t_1,\,t_2,\,t_3\in [0,T].
\end{equation}

The main result of this paper is the following uniqueness result:

\begin{thm}\label{osgood}
If the vectorfield $V$ satisfies (O) and (B),
 then there is uniqueness for
(PDE) in the class of bounded signed measures. More precisely, if
$\mu_t$ is a solution of (PDE) such that $|\mu_t|(\Rm^d)\in
L^\infty(0,T)$ then 
\begin{equation}\label{representation}
\mu_t=X(0,t,\cdot)_\#\mu_0\quad \text{ for all }\quad t\in (0,T).
\end{equation}
\end{thm}

 In the particular case when $V_t$ is continuous,
(\ref{representation})  defines a solution of (PDE) with initial
condition $\mu_0$, so that
Theorem~\ref{osgood} can also be read as an \emph{existence} result.

The same proof would give uniqueness in the larger class of measures
$\mu_t$ satisfying $|\mu_t|(\Rm^d)\in L^1(0,T)$ if conditions (O)
and (B) are
given in a stronger form with $C\in L^\infty(0,T)$, we leave the
(easy) details to the reader.

If $\rho(s)=s$, the result is well-known. It has been proved by
Bahouri and Chemin in \cite{BC} in the case $\rho(s)=s(1-\ln(s))$
(see also \cite{maniglia} for related results), under the additional
assumption that $V$ has zero divergence. The proof in \cite{BC} uses
Fourier analysis and Littlewood-Paley decompositions, and it is not
clear to us whether it can be adapted to our more general statement.

It might be tempting to think that uniqueness
for (PDE) holds in the presence of a flow of homeomorphisms
solving (ODE), but we do not know whether such a result
is true without an explicit bound like (O).

Let us now return to the definition
of the flow $X$  associated to $V$. Since $V$ is possibly discontinuous,
we consider its Filippov regularization (actually a multivalued function), namely
$$
\Vm(t,x):=\bigcap_{r>0}\overline{{\rm co}}\left(\{V(t,y) : \|y-x\|<r\}\right),
$$
where $\overline{{\rm co}}$ denotes closed convex hull. 
By definition, a generalized solution of (ODE) in the sense
of Filippov is an absolutely continuous
curve $X(t)$  such that the 
\emph{inclusion} $\dot X(t)\in\Vm(t,X(t))$ holds for almost every $t$. 
Since $x\mapsto\Vm(t,x)$ is upper semicontinuous (i.e. $x_n\to x$,
$v_n\in\Vm(t,x_n)$ and $v_n\to v$ imply $v\in\Vm(t,x)$), and
$\Vm(t,x)\neq\emptyset$, by Filippov's theorem, see \cite{Fi} or
\cite[Theorem 1.4.1]{Ho}, for all $t_1\in [0,T]$ and $x\in\Rm^d$ 
there exists a Filippov solution $X(t):[0,T]\lto \Rm^d$ 
satisfying $X(t_1)=x$.
Furthermore, $\Vm$ inherits (O) in the form
\begin{equation}\label{O1}
|\langle v-w,x-y\rangle|\leq C(t)\|x-y\|\rho(\|x-y\|)\quad
\forall v\in\Vm(t,x),\,\,\forall w\in\Vm(t,y)
\end{equation}
and this can be used to show, by the standard argument, that $X$ is unique.
The flow $X(s,t,x)$ can now be defined by requiring that
$t\mapsto X(s,t,x)$ is the only Filippov solution which
satisfies $X(s)=x$.

The following strong form of uniqueness is essential for the 
proof :
\begin{lem}\label{strong}
Let $V(t,x)$ be a vector-field satisfying (O) and (B).
Let $\gamma(s)=(t(s),x(s)):[0,L]\lto [0,T]\times \Rm^{d}$
be a Lipschitz curve such that
$$
\dot x(s) = \dot t(s)
V(t(s),x(s))
$$
for almost every $s$. If $\int_0^L |\dot t(s)| C(t(s))ds<\infty$,
where $C(t)$ is the function appearing in (O),
then 
$$
x(s)=X(t(0),t(s),x(0)).
$$ 
\end{lem}

\proof
Let us notice that the curve
$y(s):=X(t(0),t(s),x(0))$ satisfies 
$\dot y(s)\in \dot t(s)\Vm(t(s),y(s))$
for almost all $s$, and therefore \eqref{O1} gives
$$
|\langle \dot y(s)-\dot t(s)V(t(s),x),y(s)-x\rangle| \leq
|\dot t(s)|C(t(s))\|y(s)-x\|\rho(\|y(s)-x\|)
$$
for all $x\in\Rm^d$. In particular, we have
$$
|\langle \dot y(s)-\dot x(s),y(s)-x(s)\rangle| \leq |\dot
t(s)|C(t(s))\|y(s)-x(s)\|\rho(\|y(s)-x(s)\|).
$$
Denoting by $d$ the quantity
$d(s)=\|x(s)-y(s)\|$, we get
(taking into account that $\dot d(s)=0$ a.e. on $\{d=0\}$)
$$
|\dot d(s)|\leq |\dot t(s)|C(t(s))\rho(d(s))\qquad\text{for a.e.
$s\in [0,L]$}.
$$
Since the function $|\dot t(s)|C(t(s))$ is integrable, and since
$d(0)=0$, we conclude that $d(s)=0$ for all $s$.
\qed

The proof of the Theorem is now based on Smirnov's decomposition of
normal currents, see \cite{S}. We expose this theory in Section~\ref{decomposition},
and then conclude the proof of Theorem~\ref{osgood} in 
Section~\ref{proof}.

\section{Decomposition of vector fields}\label{decomposition}

Let us consider the metric space $\mL:=\text{Lip}([0,1];\Rm^k)$ of
Lipschitz curves $\gamma:[0,1]\lto \Rm^k$, endowed with the uniform
distance and the associated Borel $\sigma$-algebra. Note that the
set $\mL$ is a Borel subset of $C([0,1];\Rm^k)$, being a countable
union of compact sets. To each curve $\gamma\in\mL$, we associate
its length $L_\gamma=\int_0^1\|\dot\gamma(s)\|ds$ and the
$\Rm^k$-valued measure $T^{\gamma}=(T^\gamma_1,\ldots,T^\gamma_k)$
on $\Rm^k$ defined by
$$
\int g dT^\gamma_i=\int_0^1 g(\gamma(s))\dot\gamma_i(s) ds \qquad
i=1,\ldots,k
$$
for each bounded Borel function $g:\Rm^k\to\Rm$. Making the
supremum among all Borel functions with $\|g\|\leq 1$ we get
\begin{equation}\label{esn5}
|T^\gamma|(\Rm^k)\leq L_\gamma.
\end{equation}
Furthermore, it is easy to check that, if $\gamma$ is simple,
equality in \eqref{esn5} holds and $|T^\gamma|$ is the image of
$\|\dot\gamma\|ds$ under $\gamma$.

Let now $T=(T_1,\ldots,T_k)\in [\mM(\Rm^k)]^k$. By polar
decomposition we can write $T=W\eta$, with $W:\Rm^k\to\Rm^k$ Borel
unit vectorfield and $\eta\in\mM^+(\Rm^k)$ ($\eta$ is the total
variation of $T$ and $W$ is the orienting vectorfield, uniquely
determined up to $\eta$-negligible sets); we also assume that
$\text{div}(T)$ is (representable by) a measure
$\theta\in\mM(\Rm^k)$, namely
$$
\int \langle W,\nabla\phi\rangle d\eta=-\int \phi d\theta\qquad
\forall\phi\in C^1_c(\Rm^k).
$$
Notice that this assumption is fulfilled by $T^\gamma$, and
$$
\text{div}(T^\gamma)=\delta_{\gamma(1)}-\delta_{\gamma(0)}.
$$

We say that a measure $\nu\in \mM^+(\mL)$ is a decomposition of
$T=W\eta$ by simple curves if:
\begin{itemize}
\item[(i)] We have
\begin{equation}\label{esn3}
T=\int_{\mL} T^{\gamma}d\nu(\gamma),
\end{equation}
which explicitly means that
$$
\int\langle W,f\rangle\,d\eta= \int_{\mL}\left(\int_0^1 \langle
f(\gamma(s)),\dot\gamma(s) \rangle ds \right)d\nu(\gamma)
$$
for each bounded Borel function $f:\Rm^k\lto\Rm^k$.
\item[(ii)] We have
\begin{equation}\label{esn1}
\eta=\int_{\mL}|T^\gamma|d\nu(\gamma)
\end{equation}
and
\begin{equation}\label{esn2}
|\theta|=\int_{\mL} \big(\delta_{\gamma(1)}+\delta_{\gamma(0)} \big)
d\nu(\gamma).
\end{equation}
\item[(iii)] $\nu$-almost every curve $\gamma(t)$ is simple.
\end{itemize}

Notice that condition \eqref{esn1} can be interpreted by saying that
no cancellation occurs in \eqref{esn3}. Analogously, by applying
\eqref{esn3} to a gradient vectorfield $f$, we get
\begin{equation}\label{esn4}
\theta=\text{div}(W\eta)= \int_{\mL} \text{div}(T^\gamma)
d\nu(\gamma)=\int_{\mL} \big(\delta_{\gamma(1)}-\delta_{\gamma(0)}
\big) d\nu(\gamma).
\end{equation} So, also \eqref{esn2} implies
that no cancellation occurs in the integrals in \eqref{esn4}.

\begin{prop}\label{cachan}
Let $\nu\in\mM^+(\mL)$ be a decomposition of $W\eta$
by simple curves. Then, for
$\nu$-a.e. curve $\gamma$, we have
$|T^\gamma|=\gamma_\#(\|\dot\gamma\|ds)$ and
\begin{equation}\label{esn6}
{\dot\gamma(s)}=\|\dot \gamma(s)\|W(\gamma(s)) \qquad\text{for
a.e. $s\in [0,1]$. }
\end{equation}
\end{prop}
{\bf Proof.}
The equality $|T^\gamma|=\gamma_\#(\|\dot\gamma\|ds)$ follows
from the fact that $\nu$-almost every curve is simple.
 Inserting $f=W$ in \eqref{esn3} and taking \eqref{esn1}
into account we get
$$
\int_{\mL} \biggl(|T^\gamma|(\Rm^k)-\int_0^1\langle
W(\gamma(s)),\dot\gamma(s)\rangle ds\biggr) d\nu(\gamma)=0.
$$
Since $\nu$-almost every curve is simple, 
we have equality in  \eqref{esn5}, and we get
$$
\int_{\mL} \biggl(L_{\gamma}-\int_0^1\langle
W(\gamma(s)),\dot\gamma(s)\rangle ds\biggr) d\nu(\gamma)=0,
$$
so that
$$
\int_{\mL} \biggl(\int_0^1\|\dot\gamma(s)\|-\langle
W(\gamma(s)),\dot\gamma(s)\rangle ds\biggr) d\nu(\gamma)=0.
$$
The integrand being nonnegative, we get (\ref{esn6}).
 \qed

We can now state Theorem C of \cite{S}:
\begin{thm}\label{TS}
Any $T=W\eta$ as above can be decomposed as $\eta=\eta^0+\tilde
\eta$, where $\text{div}(W\eta^0)=0$ and $W\tilde\eta$ admits a
decomposition  $\nu\in\mM^+(\mL)$ by simple  curves.
\end{thm}

It turns out that also the divergence-free part $W\eta^0$ admits a
decomposition in ``elementary'' vector fields $T^\gamma$, but the
underlying curves $\gamma$ need not be in $\mL$: in order to obtain
the decomposition, also curves associated to Bohr quasiperiodic maps
$\gamma:\Rm\to\Rm^k$ should be considered, see \cite{S} for a
precise discussion.

\section{Proof of Theorem \ref{osgood}}\label{proof}
Let $\mu_t$ be a solution of (PDE) with initial condition
$\mu_0$ and let $S\in (0,T]$.
We want to prove that $\mu_S=X(0,S,.)_{\sharp}\mu_0$.

Let $\sigma(t,x):(0,T)\times \Rm^d\lto \{-1,1\}$ be the sign of
$\mu_t$. By this we mean a Borel map such that $\sigma\mu=|\mu|$.
Note that we really consider here a functions $\sigma$
defined at each point, and not only a class of functions
up to $|\mu|$-almost everywhere equality.
There is not a unique choice for the function $\sigma$,
but we pick one once and for all.
Let us define the vectorfield
$$
W(t,x)= \frac{\sigma(t,x)}{\|(1,V(t,x))\|}(1,V(t,x))
$$
and the Borel non-negative measure
$$\eta(t,x)=\chi_{(0,S)\times\Rm^d}\|(1,V(t,x))\|(dt\otimes |\mu_t|)
$$
on $\Rm^{d+1}=\Rm\times\Rm^d$. Note that (PDE) with initial
condition $\mu_0$ at $t=0$ and final condition $\mu_S$ at $t=S$ can
be rephrased as
$
\text{div}(W\eta)=\theta
$
in the sense of distributions in $\Rm^{d+1}$, where 
$\theta=
\delta_S \otimes \mu_S-\delta_0\otimes \mu_0$.

Let now $\eta=\eta_0+\tilde\eta$ be the decomposition provided by
Theorem~\ref{TS}, and let $\nu\in\mM^+(\mL)$ be a decomposition of
$W\tilde\eta$, with $k=1+d$. By Proposition~\ref{cachan}, $\nu$-a.e. curve
$\gamma=(t,x)$ satisfies 
$$
\dot  t(s) =\|\dot\gamma(s)\|
\frac{\sigma(t(s),x(s))}{\|(1,V(t(s),x(s)))\|}\quad,
\quad
\dot x(s)=\|\dot\gamma(s)\|\frac{\sigma(t(s),x(s))}{\|(1,V(t(s),x(s)))\|}
V(t(s),x(s)).
$$
Let us prove that for $\nu$-a.e. curve $\gamma=(t,x)$  the
integrability property
$$\int_0^1 |\dot t(s)||C(t(s))|ds<\infty
$$
holds. Indeed, take $f(t,x)=C(t)/\|(1,V(t,x))\|$ and observe that
\eqref{esn1} gives
\begin{eqnarray*}
\int_{\mL} \int_0^1|\dot t(s)||C(t(s))|ds d\nu(\gamma) &=&
 \int_{\mL}\int f d|T^\gamma| d\nu(\gamma)=\int
fd\tilde \eta\\&\leq& \int fd\eta= \int_0^T
C(t)|\mu_t|(\Rm^d)dt<\infty .
\end{eqnarray*}
In view of Lemma \ref{strong}, we conclude that
$\nu$-almost every curve
$(t,x)$ in $\mL$
satisfies
\begin{equation}\label{explicit}
x(s)=X(t(0),t(s),x(0)).
\end{equation}
Since $\nu$-almost every curve is one to one, 
we conclude that $t(s)$ is one to one for $\nu$-almost every curve.
By \eqref{esn2} we know that
$t(0)\in \{0,S\}$ and $t(1)\in \{0,S\}$ for $\nu$-almost every curve 
$\gamma=(t,x)$, and therefore, either $t(0)=0$
and $t(1)=S$, or $t(0)=S$ and $t(1)=0$.

Denoting by $\mL^+$ the Borel subset of $\mL$
formed by curves $\gamma=(t,x)$ such that 
$t$ is increasing on $[0,1]$ and satisfies $t(0)=0$
and $t(1)=S$, and by $\mL^-$
the Borel subset of $\mL$
formed by curves $\gamma=(t,x)$ such that 
$t$ is decreasing on $[0,1]$ and satisfies $t(0)=S$
and $t(1)=0$,
we conclude that 
$\nu(\mL^+\cup\mL^-)=1$.
We  denote by $\nu^\pm$ the restrictions of $\nu$
to $\mL^\pm$.
The measures $\nu^\pm$ are mutually singular, non-negative,
and $\nu=\nu^++\nu^-$.
Let 
$$B_i:\mL^+\cup \mL^-\lto \Rm^d
$$
be the Borel map defined by
$B_i(\gamma)=x(0)$ if $\gamma\in \mL^+$
and 
$B_i(\gamma)=x(1)$ if $\gamma\in \mL^-$.
Similarly, we define 
$B_f:\mL^+\cup \mL^-\lto \Rm^d
$
by
$B_i(\gamma)=x(0)$ if $\gamma\in \mL^-$
and 
$B_i(\gamma)=x(1)$ if $\gamma\in \mL^+$.
Note that 
$$B_f=B_i\circ X(0,S,\cdot)$$
$\nu$-almost everywhere by (\ref{explicit}).
Since $\theta=\delta_S \otimes \mu_S-\delta_0\otimes \mu_0$,
it follows from (\ref{esn4})
that $\mu_0=(B_i)_{\sharp}(\nu^--\nu^+)$
and  $\mu_S=(B_f)_{\sharp}(\nu^--\nu^+)$.
As a consequence, we have 
$$\mu_S=X(0,S,\cdot)_{\sharp}\mu_0.$$

\begin{small}

\end{small}
\end{document}